\documentclass[conference,letterpaper]{IEEEtran}
\IEEEoverridecommandlockouts

\usepackage{amsmath}
\usepackage[T1]{fontenc}
\usepackage[latin9]{inputenc}
\usepackage[english]{babel}
\usepackage{verbatim}
\usepackage{mathrsfs}
\usepackage{amstext}
\usepackage{amsthm}
\usepackage{amssymb}
\usepackage{esint}
\usepackage{enumerate}

\theoremstyle{plain}
\newtheorem{thm}{\protect\theoremname}
\theoremstyle{remark}
\newtheorem{rem}[thm]{\protect\remarkname}
\theoremstyle{definition}
\newtheorem{defn}[thm]{\protect\definitionname}
\theoremstyle{plain}
\newtheorem{prop}[thm]{\protect\propositionname}
\theoremstyle{plain}

\theoremstyle{plain}

\theoremstyle{plain}

\theoremstyle{definition}

\theoremstyle{definition}

\numberwithin{thm}{section}

\usepackage{tikz}
\usetikzlibrary{shapes,arrows}

\newcommand{\eps}{\varepsilon}
\newcommand{\R}{\mathbb{R}}
\newcommand{\N}{\mathbb{N}}
\newcommand{\Z}{\mathbb{Z}}
\newcommand{\mB}{\mathcal{B}}\newcommand{\mE}{\mathcal{E}}\newcommand{\mS}{\mathcal{S}}

\renewcommand{\r}{\mathrm{r}}

\newcommand{\E}{\mathbb{E}}
\newcommand{\mX}{\mathcal{X}}
\newcommand{\mP}{\mathcal{P}}
\newcommand{\bp}{\begin{proof}}
	\newcommand{\ep}{\end{proof}}
\newcommand{\ummdim}{\overline{\mathrm{mdim}}_M}

\newcommand{\udim}{\overline{\mathrm{dim}}_B}

\newcommand{\mbdim}{\overline{\mathrm{mdim}}_B}

\newcommand{\uid}{\overline{\mathrm{ID}}}

\newcommand{\mdim}{\mathrm{mdim}}

\newcommand{\id}{\mathrm{id}}

\newcommand{\LIN}{\mathrm{LIN}}

\providecommand{\conjecturename}{Conjecture}
\providecommand{\corollaryname}{Corollary}
\providecommand{\definitionname}{Definition}
\providecommand{\examplename}{Example}
\providecommand{\lemmaname}{Lemma}
\providecommand{\problemname}{Problem}
\providecommand{\propositionname}{Proposition}
\providecommand{\remarkname}{Remark}
\providecommand{\theoremname}{Theorem}

\ifCLASSINFOpdf

\fi

\begin{document}

\title{New Uniform Bounds for Almost Lossless Analog Compression}

\author{\IEEEauthorblockN{Yonatan Gutman\IEEEauthorrefmark{1}, Adam \'Spiewak\IEEEauthorrefmark{2} }
\IEEEauthorblockA{\IEEEauthorrefmark{1} Institute of Mathematics, Polish Academy of Sciences,
	ul. \'Sniadeckich 8, 00-656 Warszawa, Poland}

\IEEEauthorblockA{\IEEEauthorrefmark{2} Institute of Mathematics, University of Warsaw,
ul. Banacha 2, 02-097 Warszawa, Poland\\
Emails: y.gutman@impan.pl, a.spiewak@mimuw.edu.pl}
\thanks{We are grateful to Amos Lapidoth, Neri Merhav and Erwin Riegler for helpful discussions. Y.G was partially supported by the National Science Center (Poland) Grant 2013/08/A/ST1/00275. Y.G and A.\'{S} were partially supported by the National Science Center (Poland) grant 2016/22/E/ST1/00448.}}

\maketitle

\begin{abstract}
Wu and Verdú developed a theory of almost lossless analog compression,
where one imposes various regularity conditions on the compressor
and the decompressor with the input signal being modelled by a (typically
infinite-entropy) stationary stochastic process. In this work we consider all stationary stochastic processes with trajectories in a prescribed set $\mS \subset [0,1]^\Z$ of (bi)infinite sequences and find uniform lower and upper bounds for certain compression rates in terms of \textit{metric mean dimension} and \textit{mean box dimension}. An essential tool is the recent Lindenstrauss-Tsukamoto variational principle expressing metric mean dimension in terms of rate-distortion functions.
\end{abstract}

\textit{A full version of this paper is accessible as \cite{GS18} (preprint).}

\IEEEpeerreviewmaketitle

\section{Introduction}\label{sec:intro}
In recent years, the theory of compression for analog sources (i.e. stochastic processes with values in $\R^\Z$) underwent a major development (as a sample of such results see \cite{WV10}, \cite{jalali2017universal}, \cite{SRAB17}, \cite{GK17}). There are two key differences with the classical Shannon's model of compression for discrete sources. The first one is the necessity to employ regularity conditions on the compressor and/or decompressor functions (e.g. Lipschitz or H\"{o}lder continuity). This requirement makes the problem non-trivial and reasonable from the point of view of applications (as it induces robustness to noise). The second difference is the fact that non-discrete sources have in general infinite Shannon entropy rate, hence a different measure of complexity for stochastic processes has to be considered. One of the most fruitful approaches taken in the literature is to assume a specific structure of the source signal - as in compressed sensing, where the input vectors are assumed to be sparse (e.g. \cite{CRT06}, \cite{Donoho06}). In this setting, the theory of linear compression with efficient and stable recovery algorithms has been developed. However, strong assumptions posed on the structure of the signal reduce the applicability of the technique. A different approach was developed in the pioneering work \cite{WV10}. Instead of making assumptions on the structure of the signal, new measures of complexity related to Minkowski (box-counting) dimension of the signal were introduced and proved to be bounds on compression rates for certain classes of compressors and decompressors. Similarly, Jalali and Poor (\cite{jalali2017universal}) developed a theory of universal compressed sensing, where the linear compression rate is given in terms of a certain generalization of the R\'{e}nyi information dimension for stochastic processes with the $\psi^*$-mixing property.\\
The goal of this paper is twofold. We adapt the setting from \cite{WV10}, but instead of a single process we consider all stationary stochastic processes with trajectories in a prescribed set $\mS \subset [0,1]^\Z$. This corresponds to an a priori knowledge of all the possible trajectories of the process rather than its distribution. We deal with the question of calculating minimal compression rates in the sense of \cite{WV10} sufficient for all such stochastic processes with Borel or linear compressors and H\"{o}lder or Lipschitz decompressors. We depart from the precise setting of \cite{WV10} in several directions. We consider processes with trajectories in $[0,1]^{\Z}$, instead of $\R^\Z$ together with compression and decompression both dependent on the distribution of the process and independent of it (but dependent on $\mS$). We also consider the case where the decompressor functions are $(L, \alpha)$-H\"{o}lder with fixed $L>0$ and $\alpha \in (0,1]$ for all block lengths. Our main results are upper and lower bounds for such rates in terms of certain geometric and dynamical characteristics of the considered set $\mS$. This constitutes the second goal of the paper: we introduce notions from the theory of dynamical systems to the study of analog compression rates. As we consider stationary processes, it is natural to assume the set $\mS$ to be invariant under the shift transformation and hence it can be considered as a topological dynamical system. The obtained lower bounds are given in terms of the metric mean dimension of the system $(\mS, \mathrm{shift})$ - a geometrical invariant of dynamical systems introduced and studied by Lindenstrauss and Weiss in \cite{LW00}. Existence of connections between signal processing and mean dimension theory was observed first in \cite{GutTsu16}, where the use of the Whittaker-Nyquist-Kotelnikov-Shannon sampling theorem was essential for proving the embedding conjecture of Lindenstrauss. Another connection between these domains was established recently in \cite{lindenstrauss_tsukamoto2017rate}, where a variational principle for metric mean dimension was given in terms of rate-distortion functions. It is our main tool in developing lower bounds on compression rates for all stationary processes supported in $\mS$. In the scenario where the compressor and decompressor functions are required to be independent of the distribution of the input process (only depending on $\mS$), we introduce mean box dimension of $\mS$ as the upper bound for corresponding compression rates.

\section{Preliminaries}\label{sec:notation}
In this paper, we apply results from the theory of dynamical systems to the theory of signal processing. In line with the signal processing perspective, we consider a stationary stochastic process $\{ X_n \}_{n \in \Z},\ X_n : \Omega \to [0,1]$ defined on some probability space $(\Omega, \mathbb{P})$. Usually, instead of a single process, we are interested in considering all the stationary processes with trajectories in some prescribed set. A natural model for the set of possible trajectories is the notion of a subshift - a certain type of dynamical system. Introducing it allows us to consider stationary processes in terms of the theory of dynamical systems.

By a (topological) \textbf{dynamical system} we understand a triple $(\mX, T, \rho)$, where $(\mX, \rho)$ is a compact metric space and $T: \mX \to \mX$ is a homeomorphism. For a (countably-additive) Borel measure $\mu$ on $\mX$, by $T_* \mu$ we denote its transport by $T$, i.e. a Borel measure on $\mX$ given by $T_*\mu(A) = \mu(T^{-1}(A))$ for Borel $A \subset \mX$. We say that measure $\mu$ is \textbf{$T$-invariant}, if $\mu = T_* \mu$. By $\mP_{T}(\mX)$ we denote the set of all $T$-invariant Borel probability measures on $\mX$. We call a measure $\mu \in \mP_{T}(\mX)$ \textbf{ergodic} if every Borel set $A \subset \mX$ satisfying $T^{-1}(A) = A$ is of either full or zero measure $\mu$. The set of all ergodic measures for a transformation $T$ is denoted by $\mE_T(\mX)$. For an introduction to topological dynamics and its connections with ergodic theory see \cite[Chapters 5-8]{W82}.

Consider the unit interval $[0,1]$ with the standard metric. By the Tychonoff's theorem, $[0,1]^{\Z}$ is a compact metrizable space when endowed with the product topology. This topology is metrizable by the metric $\tau(x,y)=\sum \limits_{i=-\infty}^{\infty}\frac{1}{2^{|i|}}|x_{i}-y_{i}|$, where $x = (x_i)_{i \in \Z},\ y = (y_i)_{i \in \Z}$. This choice of the metric may seem arbitrary, but it turns out that the metric mean dimension for subshifts takes a natural form when calculated with respect to $\tau$ (see Proposition \ref{prop:canonical}). Define the shift transformation $\sigma: [0,1]^\Z \to [0,1]^\Z$ as $\sigma((x_{i})_{i=-\infty}^{\infty})=(x_{i+1})_{i=-\infty}^{\infty}.$
We are interested in properties of a given \textbf{subshift}, i.e. a closed (in the product topology) and shift-invariant subset $\mS \subset [0,1]^{\Z}$, which we interpret as the set of all admissible trajectories that can occur as input. Note that there is a one-to-one correspondence between measures $\mu \in \mP_{\sigma}(\mS)$ and distributions of stationary processes such that $(X_n)_{n \in \Z}$ belongs to $\mS$ with $\mu$-probability one. Our goal is to relate compression properties of measures (stationary processes) from $\mP_{\sigma}(\mS)$ to the geometrical properties of the set $\mS$. For $n \in \N$ define the projection $\pi_{n}:\mS\to [0,1]^{n}$ as $ \pi_n(x) = x|_{0}^{n-1} = (x_0, x_1, ..., x_{n-1}). $
For vectors $x, y \in [0,1]^n$ and $p \in [1, \infty)$, define the (normalized) \textbf{$\ell^p$ distance} as $ \| x - y\|_p=\Big(\frac{1}{n}\sum_{k=0}^{n-1}|x_{k}-y_{k}|^p\Big)^{\frac{1}{p}}$ and $\| x - y\|_\infty = \max \{ |x_k - y_k| : 1 \leq k \leq n \}.$

\section{Mean dimensions}\label{sec:dim}
In this section we will define metric mean dimension (for general dynamical systems) and (measurable) mean box dimension (for subshifts of $[0,1]^\Z$). These notions attempt to capture the average number of dimensions per iterate required to code orbits of the system. They serve as complexity measures employed to bound certain compression rates of subshifts in $[0,1]^\Z$. Let us begin with the non-dynamical notion of box dimension.
\begin{defn}\label{def:hashtag} Let $(\mX, \rho)$ be a compact metric space. For $\eps > 0$, the \textbf{$\eps$-covering number} of a subset $A \subset \mX$, denoted by $\#(A,\rho,\varepsilon)$, is the minimal cardinality $N$ of an open cover $\{U_{1},\dots,U_{N}\}$ of $A$ by sets with diameter smaller than $\varepsilon$.
\end{defn}

\begin{defn}
	Let $(\mX, \rho)$ be a compact metric space. The upper \textbf{box (Minkowski) dimension} of $A\subset\mX$
	is defined as
	\[\udim(A)=\limsup \limits_{\varepsilon\to0}\frac{\log\#(A,\rho,\varepsilon)}{\log\frac{1}{\varepsilon}}.\]
\end{defn}

In the sequel we consider only sets $A \subset [0,1]^n$ with distance induced by the norm $\| \cdot \|_{\infty}$. For more on box dimension see \cite{falconer2004fractal} and \cite{Rob11}.

\begin{defn}
	\label{def:d_n}Let $(\mX, \rho)$ be a compact metric space and let $T:\mathcal{X}\to \mathcal{X}$
	be a homeomorphism. For $n\in\N$ define a metric $\rho_{n}$ on $\mathcal{X}$
	by $\rho_{n}(x,y)=\max \limits_{0\leq k<n}\rho(T^{k}x,T^{k}y)$. Set:
	\[
	S(\mathcal{X},T,\rho,\eps)=\lim\limits_{n\to\infty}\frac{\log\#(\mathcal{X},\rho_{n},\eps)}{n}
	\]
	(the limit exists due to the subadditivity of the function $n\mapsto\log\#(\mathcal{X},\rho_{n},\eps)$).
\end{defn}
\begin{defn}\label{d:mmdim}
	The upper \textbf{metric mean dimensions} of the system $(\mX, T, \rho)$ is defined as
	\[
	\overline{\mdim}_{M}(\mathcal{X},T,\rho)=\limsup_{\eps\to0}\frac{S(\mathcal{X},T,\rho,\eps)}{\log\frac{1}{\eps}}.
	\]
\end{defn}
\begin{rem}
	It is easy to see that any system of finite topological entropy (see \cite[Chapter 7]{W82}) satisfies $\overline{\mdim}_{M}(\mathcal{X},T,\rho)=0$. Metric mean dimension can be easily computed for full shifts: if $(A,d)$ is a compact metric space, then $\ummdim(A^\Z, \sigma, \rho) = \udim(A, d)$, where $\rho$ is the product metric (see \cite{GS18}). Also, $\ummdim$ is an invariant for bi-Lipshitz isomorphisms: if $(\mX, T, \rho_1)$ and $(\mathcal{Y}, S, \rho_2)$ are dynamical systems and $\Phi : \mX \to \mathcal{Y}$ is bi-Lipshitz and equivariant (i.e. $\Phi \circ T = S \circ \Phi$), then $\ummdim(\mX, T, \rho_1) = \ummdim(\mathcal{Y}, S, \rho_2)$.
\end{rem}

A topological version of mean dimension for actions of amenable groups was introduced by Gromov in \cite{G} and studied by Lindenstrauss and Weiss in their seminal work \cite{LW00}. It turns out that the topological mean dimension is the right invariant to study for the problem of existence of an embedding into $(([0,1]^D)^\Z, \sigma)$ (see \cite{GutTsu16}). For more on mean topological dimension see \cite{coornaert2015topological}. The metric mean dimension was introduced in \cite{LW00} and proved to be, when calculated with respect to any compatible metric, an upper bound for the topological mean dimension.

When $\mS\subset[0,1]{}^{\mathbb{{Z}}}$ is a subshift and $\rho=\tau$ (see
Section \ref{sec:notation}), metric mean dimension can be expressed in a more canonical
form:
\begin{prop}\label{prop:canonical} For a subshift $\mS\subset[0,1]^{\Z}$ it holds
	\[
	\overline{\mdim}_{M}(\mS,\sigma,\tau)=\limsup_{\varepsilon\to0}\lim_{n\rightarrow\infty}\frac{\log\#(\pi_{n}(\mS),||\cdot||_{\infty},\varepsilon)}{n\log\frac{1}{\varepsilon}}.
	\]
\end{prop}

\begin{defn}
For $\mS\subset[0,1]{}^{\mathbb{{Z}}}$ we define its upper \textbf{mean box dimension} as
\[
\mbdim(\mS)=\lim_{n\to\infty}\frac{\overline{\dim}_{B}(\pi_{n}(\mS))}{n},
\]
where $\udim(\pi_n(\mS))$ is calculated with respect to $\| \cdot \|_\infty$ norm on $[0,1]^n$. The limit exists due to the subadditivity of the function $n\mapsto\overline{\dim}_{B}(\pi_{n}(\mS))$.
\end{defn}
\begin{prop}
	\label{prop:mdim leq mbdim} Let $\mS\subset[0,1]^{\Z}$ be a subshift. Then
	\[
	\overline{\mdim}_{M}(\mS, \sigma,\tau)\leq\mbdim(\mS).
	\]
\end{prop}

In \cite{WV10}, Wu and Verdú gave bounds on certain compression rates in terms of the following notion.
\begin{defn}\label{d:measurable_mean_minkowski}\cite[Def. 10]{WV10}) For a subshift $\mS\subset[0,1]^{\Z}$, invariant measure $\mu\in \mP_{\sigma}(\mS)$,
	$n\in\N$ and $0 \leq \delta<1$ define the \textbf{measurable mean box dimension} as
	\begin{multline*}
	R_{B}(\mu,\delta)= \limsup \limits_{n \to \infty}\ \inf \Big\{ \frac{\overline{\dim}_{B}(A)}{n} :  A\subset[0,1]^{n},\\
	A \text{ - compact, } \mu(\pi_{n}^{-1}(A))\geq 1-\delta\Big\}.
	\end{multline*}
\end{defn}

\begin{rem}
Wu and Verdú use the name \textit{Minkowski-dimension compression rate} for $R_{B}(\mu,\delta)$ . As we reserve the term \textit{compression rate} for a different concept (of an operational meaning, see Section \ref{sec:compression_rates}), we decided to introduce a different name.
\end{rem}
\section{Analog compression}\label{sec:main}

In this section we introduce analog compression rates for sources with alphabet $[0,1]$ and state our main results. In this setting it is natural to assume regularity constraints on the compressor and decompressor functions. This follows from the fact that we are taking an infinite alphabet under consideration: for every $n \in \N$ there exists a (Borel) bijection between $[0,1]^n$ and $[0,1]$, hence the corresponding compression rates tend to zero if we do not assume any further regularity of the compressor and decompressor functions (cf. \cite[Section IV.B]{WV10}). On the other hand, from the point of view of applications it is desirable to impose some regularity conditions, as they induce robustness to noise and enable numerical control of the errors occurring in the compression and decompression processes.

\subsection{Compression rates}\label{sec:compression_rates}
\begin{defn}
	\label{def:(regularity-classes)}
	A regularity class is a set $\mathcal{C}$ of functions between finite dimensional unit cubes, i.e. $\mathcal{C} \subset \{ f : [0,1]^n \to [0,1]^k\ |\ n,k \in \N \}$.
\end{defn}
We will consider the following regularity classes: 	$\mathcal{B} = \{ \text{Borel maps}\}$, $\mathcal{{H}}_{\alpha} = \{ \alpha \text{-H\"{o}lder maps}\}$, $\mathcal{{H}}_{L,\alpha} = \{ \alpha \text{-H\"{o}lder maps with constant } L \}$, $\mathrm{LIN} =  \{ \text{linear maps} \},$ where the H\"{o}lder condition is considered with respect to $\|\cdot\|_{\infty}$ on $[0,1]^n$ and $[0,1]^k$. Below we define several compression rates for various requirements on the performance of the compression and decompression process (see also \cite[Def. 3]{WV10}).

\begin{defn}
	\label{def:ratio is achievable} Let
	$\mS \subset [0,1]^{\mathbb{Z}}$ be a subshift
	and $\mu\in  \mP_{\sigma}(\mS)$. Let $\mathcal{{C}},\mathcal{{D}} \subset \{ f : [0,1]^n \to [0,1]^k\ |\  n, k \in \N \}$
	be regularity classes. For $n \in \N$ and $\eps \geq 0$, the $\mathcal{C}-\mathcal{D}$ \textbf{almost lossless analog compression rate $\r_{\mathcal{C}-\mathcal{D}}(\mu,\eps,n)\geq0$} of $\mu$ with $n$-block error probability $\eps$ is the infimum of $\frac{k}{n}$,	where $k$ runs over all natural numbers such that there exist maps $f:[0,1]^{n}\rightarrow[0,1]^{k},\ f\in\mathcal{{C}}$ and $g:[0,1]^{k}\rightarrow[0,1]^{n},\ g\in\mathcal{{D}}$
	with
	\begin{equation}\label{eq:compression_rate_def}
	\mu(\{x\in \mS|\ g\circ f(x|_{0}^{n-1})\neq x|_{0}^{n-1}\})\leq\eps.
	\end{equation}
	Define further $
	\r_{\mathcal{C}-\mathcal{D}}(\mu,\eps)=\limsup \limits_{n\rightarrow\infty}\ \r_{\mathcal{C}-\mathcal{D}}(\mu,\eps,n).$
\end{defn}

We define similarly the $\mathcal{C}-\mathcal{D}$ \textbf{uniform almost lossless analog compression rate $\r_{\mathcal{C}-\mathcal{D}}(\mS,\eps)\geq0$} of $\mS$ by requiring that (\ref{eq:compression_rate_def}) holds \textit{for all} $\mu \in \mP_{\sigma}(\mS)$. In such a case, compression can be performed at asymptotic rate $\r_{\mathcal{C}-\mathcal{D}}(\mS, \eps)$ without knowing the distribution from which data comes, as long as the process is supported in $\mS$.

For $p\geq 1$ we define also the $\mathcal{C}-\mathcal{D}$ \textbf{probability analog compression rate $\r_{\mathcal{C}-\mathcal{D}}^{P,p}(\mu,\eps,n, \delta)\geq0$} of $\mu$ with $n$-block error probability $\delta \geq 0$ at scale $\eps$ by replacing condition (\ref{eq:compression_rate_def}) with
\begin{equation}\label{e:compression_rate_condition_prob}
\mu( \{ x \in \mS : \|x|_0^{n-1} - g \circ f (x|_0^{n-1})\|_p \geq \eps \}) \leq \delta.
\end{equation}
We define further $
\r_{\mathcal{C}-\mathcal{D}}^{P, p}(\mu,\eps, n)=\lim \limits_{\delta \to 0}\ \r_{\mathcal{C}-\mathcal{D}}^{P,p}(\mu,\eps,n, \delta)$ and $\r_{\mathcal{C}-\mathcal{D}}^{P,p}(\mu,\eps)=\limsup \limits_{n\rightarrow\infty}\ \r_{\mathcal{C}-\mathcal{D}}^{P,p}(\mu,\eps,n)$. We do not use $\r_{\mathcal{C}-\mathcal{D}}^{P,p}$ directly in this paper, but it allows us to state results of \cite{jalali2017universal} in the language of compression rates.

\subsection{Previous results}
Let us begin by presenting some known results giving bounds on compression rates introduced in the previous subsection. In their pioneering article \cite{WV10} Wu and Verdú calculated and gave bounds on $\r_{\mathcal{{C}}-\mathcal{{D}}}(\mu,\eps)$ for
certain $\mathcal{{C}}$ and $\mathcal{{D}}$ and fixed $\mu\in\mathcal{{P}}_{\sigma}(\mathbb{{R}^{\mathbb{{N}}}})$.
For example by \cite[Thm.  9]{WV10} it follows for Bernoulli measure
$\mu=\bigotimes\limits_{\mathbb{{Z}}}\nu\in\mathcal{{P}}_{\sigma}([0,1]{}^{\mathbb{{Z}}})$ that $\r_{\mathcal{{B}}-\mathcal{H}_1}(\mu,\eps)\ge \uid(\nu)$ for $0<\eps<1$, where $\uid$ denotes the upper Rényi information dimension of a probability measure. Another of their results is the following:

\begin{thm}\label{thm:wu_verdu_up}\cite[Thm.  18]{WV10} For $\mu \in \mP_{\sigma}([0,1]^{\Z})$ and $\alpha \in (0,1)$ the following holds: \[\r_{\mathrm{LIN}-\mathcal{H}_{\alpha}}(\mu, \eps) \leq \frac{1}{1 - \alpha}R_B(\mu, \eps)\]
and consequently $\r_{\mathrm{LIN}-\mathcal{H}}(\mu, \eps) \leq R_B(\mu, \eps)$.
\end{thm}

\begin{rem}
	The above upper bound on $\r_{\mathrm{LIN}-\mathcal{H}_{\alpha}}(\mu, \eps)$ comes from minimizing $R$ in \cite[(172)]{WV10} for fixed $\beta$. Stronger result than the existence of linear compressor and H\"{o}lder decompressor was proven in \cite[Section VIII]{SRAB17}, where it is shown that almost every linear transformation of rank large enough serves as a good compressor in this setting.
\end{rem}

For the other direction, following closely the proof of the upper bound in \cite[Equation (75)]{WV10}, we have the following proposition (see \cite{GS18} for the proof).
\begin{prop}\label{prop:rbh_below}
	Let $\mS \subset [0,1]^\Z$ be a subshift
	and $\mu\in  \mP_{\sigma}(\mS)$. Then $
	\alpha R_{B}(\mu,\delta)\leq \r_{\mathcal{B}-\mathcal{H_{\alpha}}}(\mu,\delta)$ for $0<\delta<1$ and $\alpha \in (0,1]$.
\end{prop}

In applications the measure governing the source is not always known. Some universality in the compression process was proposed in \cite{jalali2017universal}. In terms of compression rates, the following bound was obtained (for the definition of $\overline{d}_0(\mu)$ see \cite[Def. 2]{jalali2017universal} and for $\psi^*$-mixing see \cite[Def. 3]{jalali2017universal}):
\begin{thm}\label{thm:jalali_poor}(\cite[Thms 7,8]{jalali2017universal})
	Let $\mu \in \mP_{\sigma}([0,1]^\Z)$ be $\psi^*$-mixing. Then
	\[ \sup \limits_{\eps > 0}\ \r^{P,2}_{\LIN - \mB}(\mu, \eps) \leq \overline{d}_0(\mu). \]
\end{thm}

\begin{rem}
	\cite{jalali2017universal} proved more than merely existence of suitable linear compressors. More precisely, they proved that for any $\eta>0$, if $(X_n)_{n \in \Z}$ is a $\psi^*$-mixing stochastic process with distribution $\mu$ and $A_n \in \R^{n \times m_n}$ are independent random matrices with entries drawn i.i.d according to $\mathcal{N}(0, 1)$ and independently from $(X_n)_{n \in \Z}$ with $\frac{m_n}{n} \geq (1+ \eta) \overline{d}_0(\mu)$, then
	\[\|X|_0^{n-1} - g_n \circ A_n (X|_0^{n-1}) \|_2 \overset{n \to \infty}{\longrightarrow} 0 \text{ in probability } \mu \otimes \nu,\]
	where $\nu$ is the distribution of $(A_n)_{n=1}^{\infty}$ and $g_n : \R^{m_n} \to \R^n$ are some explicitly defined Borel functions  (depending only on $A_n$). Hence, for such a random sequence of matrices, the expected value
	\[ \E_{\nu} \mu( \{ x \in [0,1]^\Z : \|x|_{0}^{n-1} - g_n \circ A_n (x|_{0}^{n-1})\|_2 \geq \eps \}) \]
	tends to zero as $n \to \infty$ for any $\psi^*$-mixing measure $\mu \in \mP_{\sigma}([0,1]^{\Z})$. Theorem \ref{thm:jalali_poor} follows from this, since for any $\delta >0$ and $n$ large enough, there exists $A \in \R^{n \times m_n}$ satisfying
	\begin{equation}\label{e:jalali_poor} \mu( \{ x \in [0,1]^\Z : \|x|_{0}^{n-1} - g_n \circ A_n (x|_{0}^{n-1})\|_2 \geq \eps \}) \leq \delta.
	\end{equation}
	The decompressors $g_n$ take only finitely many values (hence are not continuous) and are defined via a certain minimization problem (which makes the decompression algorithm implementable, though not efficient (cf. \cite[Remark 3]{jalali2017universal})). The authors proved also that, in a certain setting, such a compression scheme is robust to noise (see \cite[Thms 9 and 10]{jalali2017universal}). The strength of the result is the universality of the compression scheme, which is designed without any prior knowledge of the distribution $\mu$: a random Gaussian matrix will serve as a good compressor as long as the rate is at least $\overline{d}_0(\mu)$. However, it does not follow that one can choose a sequence of matrices $A_n$ satisfying (\ref{e:jalali_poor}) \textit{for all} $\psi^*$-mixing measures $\mu$ with $\overline{d}_0(\mu) \leq d$ for some $d \in [0,1]$. Also, $\psi^*$-mixing is quite a restrictive assumption.
\end{rem}

\subsection{Main results}

Instead of assuming specific properties of the measure governing the source, we consider the scenario in which the set of all possibles trajectories is known. Therefore we are interested in the following question:\\ \ \\
\textbf{Main Question:} Given a subshift $\mS \subset [0,1]^\Z$, calculate
\[\sup \limits_{\eps > 0} \sup \limits_{\mu\in\mathcal{P}_{\sigma}(\mS)}\r_{\mathcal{{C}}-\mathcal{{D}}}(\mu,\eps) \text{ and } \sup \limits_{\eps > 0}\ \r_{\mathcal{{C}}-\mathcal{{D}}}(\mS,\eps)\] for fixed regularity classes $\mathcal{C}$ and $\mathcal{D}$.

We are interested in this question for $\mathcal{C} \in \{\mathcal{B}, \mathrm{LIN}\}$ and $\mathcal{D}=\mathcal{H}_{L,\alpha}$. Such or similar regularity conditions have appeared previously in the literature (e.g. Theorems \ref{thm:wu_verdu_up} and \ref{thm:jalali_poor}). As above quantities are decreasing with $\eps$, one can exchange $\sup \limits_{\eps > 0}$ for $\lim \limits_{\eps \to 0}$. Taking supremum over invariant measures in Theorem \ref{thm:wu_verdu_up} and Proposition \ref{prop:rbh_below}, we obtain:

\begin{thm}\label{thm:main_holder_vw} Let $\mS \subset [0,1]^\Z$ be a subshift. The following holds for
	every $0<\alpha<1$: \[
	\alpha \sup \limits_{\eps>0}\ \sup \limits_{\mu\in  \mP_{\sigma}(\mS)} R_B (\mu, \eps)\leq \sup \limits_{\eps>0}\ \sup \limits_{\mu\in  \mP_{\sigma}(\mS)}\r_{\mathcal{{LIN}}-\mathcal{{H}}_{\alpha}}(\mu,\eps) \leq\]
	\[  \leq \frac{1}{1-\alpha}\sup \limits_{\eps>0}\ \sup \limits_{\mu\in  \mP_{\sigma}(\mS)} R_B (\mu, \eps).
	\]
\end{thm}

Note that the above results do not give an explicit bound on the constant $L$; in fact, they do not guarantee a uniform bound for $L$ among the sequence of decoders. This is a drawback from the point of view of error control. Hence, it is reasonable to consider also class $\mathcal{H}_{L, \alpha}$ for fixed $L, \alpha$. Note that $r_{\mathcal{C} - \mathcal{H}_{\alpha}}(\mu, \eps) \leq r_{\mathcal{C} - \mathcal{H}_{L, \alpha}}(\mu, \eps)$ for any compression rate and class $\mathcal{C}$. In the sequel we give both lower and upper bounds for $\r_{\mathcal{B} - \mathcal{H}_{L, \alpha}}(\mu, \eps)$ and $\r_{\mathrm{LIN} - \mathcal{H}_{L, \alpha}}(\mu, \eps)$ in terms of $\ummdim(\mS, \sigma, \tau)$ and $\mbdim(\mS)$. Note that the quantities $R_B(\mu, \eps)$ depending on the measure and parameter $\eps$ might be harder to calculate in specific examples than various geometric mean dimensions. Our main results are the following:

\begin{thm}\label{thm:main_holder_low} Let $\mS \subset [0,1]^\Z$ be a subshift. The following holds for
	every $0<\alpha \leq 1, L>0$:
	\begin{equation*}
	\alpha \ummdim(\mS, \sigma,\tau) \leq \sup \limits_{\eps>0}\ \sup \limits_{\mu\in  \mP_{\sigma}(\mS)}\ \r_{\mathcal{{B}}-\mathcal{{H}}_{L,\alpha}}(\mu,\eps).
	\end{equation*}
\end{thm}

For a sketch of the proof see Section \ref{sec:lower}. For details and extension to $L^p$ compression rates see \cite{GS18}. In general, equality does not hold in Theorem \ref{thm:main_holder_low}. We also cannot change the class $\mathcal{H}_{L, \alpha}$ to $\mathcal{H}_{\alpha}$, i.e. $\alpha \ummdim(\mS, \sigma, \tau)$ cannot serve as a lower bound in Theorem \ref{thm:main_holder_vw}. See \cite{GS18} for suitable examples. 

\begin{thm}\label{thm:main_holder_up}
	Let $\mS \subset [0,1]^\Z$ be a subshift. Then, for
	every $0<\alpha < 1$
	\[ \sup \limits_{\eps>0}\ \sup \limits_{\mu\in  \mP_{\sigma}(\mS)}\ \inf \limits_{L>0}\ \r_{\mathrm{LIN}-\mathcal{{H}}_{L,\alpha}}(\mu,\eps) \leq \]
	\[ \leq \inf \limits_{L>0}\ \r_{\mathrm{LIN}-\mathcal{{H}}_{L,\alpha}}(\mS,0) \leq  \min \{1, \frac{2}{1-\alpha}\mbdim(\mS)\}. \]
\end{thm}
The proof is based on the embedding theorem for $\udim$ with H\"{o}lder inverse \cite[Thm.  4.3]{Rob11} (see \cite{BGS18} for an almost sure embedding theorem for Hausdorff dimension). See \cite{GS18} for the proof and examples showing that one cannot change the constant $\frac{2}{1 - \alpha}$ to $\frac{t}{1 - \alpha}$ for $t < 2$ and $\inf \limits_{L>0}$ cannot be omitted.

\section{Rate-distortion functions and variational principles for metric mean dimension}\label{sec:rate_dist}

Our proof of the lower bound in Theorem \ref{thm:main_holder_low} is based on a variational principle for metric mean dimension in terms of rate-distortion function \cite{lindenstrauss_tsukamoto2017rate}. We work with a slight modification of the expression used in \cite{lindenstrauss_tsukamoto2017rate}.

\begin{defn}\label{def:rate_distortion_function}
	(compare with \cite[p. 3-4]{lindenstrauss_tsukamoto2017rate}) Let $(A,d)$ be a compact metric space, let
	$\mS \subset A^{\mathbb{Z}}$ be a subshift and $\mu\in  \mP_{\sigma}(\mS)$. For $\varepsilon>0$ and $n\in\mathbb{{N}}$
	we define the \textbf{rate-distortion function} $\tilde{R}_{\mu}(n,\varepsilon)$
	as the infimum of $\frac{I(X;Y)}{n}$, where $X=(X_0, ..., X_{n-1})$ and $Y=(Y_{0},\dots,Y_{n-1})$ are random variables defined on some probability space $(\Omega,\mathbb{P})$ such that
	\begin{itemize}
		\item $X=(X_0, ..., X_{n-1})$ takes values in $A^n$, and
		its law is given by $(\pi_n)_*\mu$.
		\item $Y=(Y_{0},\dots,Y_{n-1})$ takes values in $A^n$ and $\mathbb{{E}}\left(\frac{1}{n}\sum_{k=0}^{n-1}d(X_{k},Y_{k})\right)\leq\varepsilon$.
	\end{itemize}
	Here $I(X;Y)$ is the mutual information of random vectors $X$ and $Y$ (see \cite{Gray11} and \cite{lindenstrauss_tsukamoto2017rate}). The function $n \mapsto n\tilde{R}_{\mu}(n, \eps)$ is subadditive (see \cite[Thm.  9.6.1]{Gallager68} for a proof in the finite alphabet case). Hence, we may define
	\[
	\tilde{R}_{\mu}(\eps)=\lim_{n \to \infty}\tilde{R}_{\mu}(n,\varepsilon) = \inf_{n \in \N}\tilde{R}_{\mu}(n,\varepsilon).
	\]
	
\end{defn}

The following theorem is a variant of the variational principle for metric mean dimension in the case of subshifts. It is deduced from the original theorem \cite[Theorem III.1]{lindenstrauss_tsukamoto2017rate}. We also prove that one can take the supremum over ergodic measures (see \cite{GS18} for details).
\begin{thm}\label{thm:var_prin_R_tilde}
	Let $\mS \subset [0,1]^\Z$ be a subshift. Then
	\begin{equation*}
	\begin{split}
	\overline{\mdim}_{\mathrm{M}}(\mS, \sigma,\tau) & =\limsup_{\varepsilon\to0}\sup \limits_{\mu\in  \mP_{\sigma}(\mS)}\frac{\tilde{R}_{\mu}(\varepsilon)}{\log\frac{1}{\eps}} = \\ & =\limsup_{\varepsilon\to0}\sup \limits_{\mu\in  \mE_{\sigma}(\mS)}\frac{\tilde{R}_{\mu}(\varepsilon)}{\log\frac{1}{\eps}}.
	\end{split}
	\end{equation*}
\end{thm}
The above theorem remains true if we consider the $L^p$ distortion function instead of the $L^1$ variant (see \cite{GS18}). As proved in \cite[Thm.  1]{GK17}, for the $L^2$ rate-distortion function the above limit \textit{for fixed} $\mu \in \mP_{\sigma}(\mS)$ gives the upper information dimension of $\mu$. For a variational principle for $\ummdim$ in terms of the mean R\'{e}nyi information dimension see \cite{GS18}.

\section{Lower bounds}\label{sec:lower}

The following inequality is the main ingredient of the proof of Theorem \ref{thm:main_holder_low}, as together with Theorem \ref{thm:var_prin_R_tilde} it yields the result. However, it is of independent interest, since it gives a lower bound for $\r_{\mathcal{{B}}-\mathcal{{H}}_{L,\alpha}}(\mu,\eps)$ for fixed $\mu$ and $\eps$.

\begin{thm}\label{thm:main}
	Let $\mS \subset [0,1]^\Z$ be a subshift.  The following
	holds for $\mu\in  \mP_{\sigma}(\mS),\ 0<\alpha\leq1,\ L>0$:
	\[
	\frac{\tilde{R}_{\mu}((\frac{L}{2^{ \alpha}}+\eps^{(1-\alpha)})\eps^{\alpha})}{\log(\lceil\frac{1}{\eps}\rceil)}\leq \r_{\mathcal{{B}}-\mathcal{{H}}_{L,\alpha}}(\mu,\eps).
	\]
\end{thm}

\begin{proof}
	Fix $\delta,\eps>0$. Assume that $\mS$
	achieves $\mathcal{{B}}-\mathcal{H}_{L, \alpha}$ almost lossless analog compression rate $\r_{\mathcal{{B}}-\mathcal{H}_{L, \alpha}}(\mu,\eps)<\infty$
	with error probability $\eps$. One may find $k,n\in\mathbb{N}$ with
	$\frac{k}{n}\leq \r_{\mathcal{{B}}-\mathcal{{H}}_{L,\alpha}}(\mu,\eps)+\delta$ and functions
	$f:[0,1]^{n}\rightarrow[0,1]^{k},\ f \in \mathcal{B}$,
	$g:[0,1]^{k}\rightarrow[0,1]^{n},\ g \in \mathcal{H}_{L, \alpha}$ such that
	$\mu(E)\leq\eps$, where $E=\{x\in \mathcal{X}|\ g\circ f(x|_{0}^{n-1})\neq x|_{0}^{n-1}\}$.
	Regularly partition $[0,1]^{k}$ into $\lceil\frac{1}{\eps}\rceil^{k}$ cubes of side
	$\lceil\frac{1}{\eps}\rceil^{-1}$ Borel-wise and let $c:[0,1]^{k}\rightarrow F$ associate to each point the	center of its cube. Note that $|F|=\lceil\frac{1}{\eps}\rceil^{k}$
	and $||x-c(x)||_{\infty}\leq\frac{\eps}{2}$ for all $x\in[0,1]^{k}$.
	Define $Y:[0,1]^{n}\rightarrow[0,1]^{n}$ by $Y(p)=g(c(f(p)))$ and $X: [0,1]^{n} \to [0,1]^{n}$ by $X = \id$. This gives a pair of random vectors on the probability space $([0,1]^n, (\pi_n)_* \mu)$.
	We now estimate (here $A=[0,1]$ and $d=\| \cdot \|_{\infty}$)
	\[
	\mathbb{{E}}\left(\frac{1}{n}\sum_{k=0}^{n-1}d(X_{k},Y_{k})\right) \leq \int \limits_{[0,1]^n} \|x - g  \circ f (x) \|_\infty d (\pi_n)_* \mu(x) + \]
	\[ + \int \limits_{[0,1]^n} \|g \circ f(x) - g \circ c \circ f (x) \|_\infty d (\pi_n)_* \mu(x) \leq \]
	\[ \leq \eps + \int \limits_{[0,1]^n} L \|f(x) - c \circ f(x)\|_{\infty}^{\alpha}d (\pi_n)_* \mu(x) \leq \eps + L \frac{\eps^{\alpha}}{2^{\alpha}}. \]
	This implies
	\[
	\tilde{R}_{\mu}((\frac{L}{2^{ \alpha}}+\eps^{1-\alpha})\eps^{\alpha})\leq\frac{1}{n}I(X;Y)\leq\frac{1}{n}H(Y)\leq\]
	\[ \leq \frac{\log(\lceil\frac{1}{\eps}\rceil^{k})}{n}=\frac{k\log(\lceil\frac{1}{\eps}\rceil)}{n}\leq \log(\lceil\frac{1}{\eps}\rceil)(\r_{\mathcal{{B}}-\mathcal{{H}}_{L, \alpha}}(\mu,\eps)+\delta).
	\]
\end{proof}

\bibliographystyle{IEEEtran}
\bibliography{IEEEabrv,universal_bib_nourl}

\end{document}